# Properties of Polytopes Representing Natural Numbers

Ya-Ping Lu and Shu-Fang Deng

East China University of Science and Technology

**Abstract**: Lattice polytope representation of natural numbers is introduced based on the fundamental theorem of arithmetic. The combinatorial and geometric properties of the polytopes are studied using *Polymake* and *Qhull* software. The volume of the polytope representing a natural number and the sum of the volumes of polytopes representing a subset of natural numbers are further examined.

Figurate numbers are a sequence of numbers associated with a certain geometric form in 1, 2, 3, or higher dimension. Polygonal and pyramidal numbers are examples of 2D and 3D figurate numbers. A figurate number presents the links between Arithmetic and Geometry [1]. The objective of this paper is to establish a link between natural numbers and geometry by representing the natural number using geometric objects. Each natural number is "assigned" to a geometric object or, more precisely, a multi-dimensional lattice polytope. The construction of the polytope is based on the fundamental theorem of arithmetic. These lattice polytopes constructed serve as the "number figures" that represent the natural numbers.

## 1. Representation of Natural Numbers by Convex Lattice Polytopes

Natural numbers are positive integers, and the set of natural numbers is denoted by

$$N = \{1, 2, 3, \dots\}.$$

The fundamental theorem of arithmetic (or factorization theorem) states that every positive integer except the number 1 can be represented uniquely as a product of one or more prime numbers [2]. For any natural number, a lattice polytope can be defined based on the factorization theorem.

**Definition 1.1. (lattice polytope representing a natural number).** Let $n$ be the number of prime numbers less than or equal to a natural number N, or $n = \pi(N)$, where $\pi(N)$ is the prime counting function. The lattice polytope representing N is defined as

$$P(N) := conv\{v_1, v_2, \ldots, v_N\}$$

where $v_M \in \mathbf{R}^n$ (M $= 1, 2, \ldots, N$) is the lattice point in the $n$-dimensional Euclidian space with the Cartesian coordinates of

$$v_M = (x_{1,M}, x_{2,M}, \ldots, x_{n,M})$$

and $M$ is expressed by the first $n$ prime numbers as

$$M = \prod_{i=1}^{n} p_i^{x_{i,M}} \quad x_{i,M} \in Z_{\geq 0}.$$

The expression of M in Definition 1.1 has exactly n prime numbers and the exponent $x_{i,M}$ is a non-negative integer, which is slightly different from the prime factorization of natural numbers normally formulated, in which all exponents are positive integers.

With Definition 1.1, every natural number $N$ is represented by a $\pi(N)$-dimensional lattice polytope constructed by $N$ lattice points associated with the $N$ consecutive natural numbers from 1 to $N$. Note that the lattice point associated with the number 1 is at the origin of the Euclidian space with coordinates of zero in all dimensions. The first 16 natural numbers and the volumes of the polytopes representing the numbers are calculated using *Qhull* (2019.1) software [3] and given in Table 1.1.

**Table 1.1 Polytopes representing the first 16 natural numbers**

| N | $n$ | Polytope, $P$ | No. of lattice points | Polytope vertices | $vol(P)$ |
|---|---|---|---|---|---|
| 1 | 0 | Point | 1 | at the origin | 1 |
| 2 | 1 | Line segment | 2 | {0, 1} | 1 |
| 3 | 2 | Triangle | 3 | {(0,0),(1,0),(0,1)} | 1/2 |
| 4 | | | 4 | {(0,0),(2,0),(0,1)} | 1 |
| 5 | 3 | Tetrahedron | 5 | {(0,0,0),(2,0,0),(0,1,0),(0,0,1)} | 1/3 |
| 6 | | Trapezoid based pyramid | 6 | {(0,0,0),(2,0,0),(0,1,0),(1,1,0),(0,0,1)} | 1/2 |
| 7 | 4 | 4-polytope | 7 | {(0,0,0,0),(2,0,0,0),(0,1,0,0),(1,1,0,0),(0,0,1,0),(0,0,0,1)} | 1/8 |
| 8 | | | 8 | {(0,0,0,0),(3,0,0,0),(0,1,0,0),(1,1,0,0),(0,0,1,0),(0,0,0,1)} | 1/6 |

| 9 | | | 9 | {(0,0,0,0),(3,0,0,0),(0,2,0,0),(0,0,1,0),(0,0,0,1)} | 1/4 |
|---|---|---|---|---|---|
| 10 | | | 10 | {(0,0,0,0),(3,0,0,0),(0,2,0,0),(0,0,1,0),(1,0,1,0),(0,0,0,1)} | 1/3 |
| 11 | 5 | 5-polytope | 11 | {(0,0,0,0,0),(3,0,0,0,0),(0,2,0,0,0),(0,0,1,0,0),(1,0,1,0,0), (0,0,0,1,0),(0,0,0,0,1)} | 1/15 |
| 12 | | | 12 | {(0,0,0,0,0),(3,0,0,0,0),(0,2,0,0,0),(2,1,0,0,0),(0,0,1,0,0), (1,0,1,0,0),(0,0,0,1,0),(0,0,0,0,1)} | 3/40 |
| 13 | | | 13 | {(0,0,0,0,0,0),(3,0,0,0,0,0),(0,2,0,0,0,0),(2,1,0,0,0,0),(0,0,1,0,0,0), (1,0,1,0,0,0),(0,0,0,1,0,0),(0,0,0,0,1,0),(0,0,0,0,0,1)} | 1/80 |
| 14 | | | 14 | {(0,0,0,0,0,0),(3,0,0,0,0,0),(0,2,0,0,0,0),(2,1,0,0,0,0),(0,0,1,0,0,0), (1,0,1,0,0,0),(0,0,0,1,0,0),(0,0,0,0,1,0),(1,0,0,1,0,0),(0,0,0,0,0,1)} | 11/720 |
| 15 | 6 | 6-polytope | 15 | {(0,0,0,0,0,0),(3,0,0,0,0,0),(0,2,0,0,0,0),(2,1,0,0,0,0),(0,0,1,0,0,0), (1,0,1,0,0,0),(0,1,1,0,0,0),(0,0,0,1,0,0),(0,0,0,0,1,0),(1,0,0,1,0,0), (0,0,0,0,0,1)} | 7/360 |
| 16 | | | 16 | {(0,0,0,0,0,0),(4,0,0,0,0,0),(0,2,0,0,0,0),(0,0,1,0,0,0),(1,0,1,0,0,0), (0,1,1,0,0,0),(0,0,0,1,0,0),(0,0,0,0,1,0),(1,0,0,1,0,0),(0,0,0,0,0,1)} | 1/45 |

## 2. Studies of Polytope Properties Using *Polymake*

The convex lattice polytopes representing the natural numbers are studied using *Polymake* [4-5] (version 3.2 and Version 4.0) with Ubuntu [6] operating system.

Every polytope representing a natural number $N > 1$ can be constructed by adding the $N$-th lattice point to the polytope representing the previous number, $N - 1$. Thus, the construction of the polytopes can be done sequentially staring from smaller numbers. Some of the properties of the polytopes representing the natural numbers up to 448 are provided in Appendix A. Polytopes for natural number bigger than 448 have not been calculated due to RAM limitation of the computer used.

The parameters of the polytopes listed in Appendix A are the following:

- Column 1: Natural number, N, which is also the number of lattice points of the polytope, N_LATTICE_POINTS (notation used in *Polymake*);
- Column 2: The diameter of the polytope, DIAMETER;
- Column 3: Dimension of the polytope, DIM;
- Column 4: Normalized volume or relative volume of the polytope, LATTICE_VOLUME;
- Column 5: The number of vertices of the polytope, N_VERTICES;
- Column 6: The number of edges of the polytope, N_EDGES;
- Column 7: The number of facets of the polytope, N_FACETS;
- Column 8: The maximal integral width of the polytope with respect to the facet normals, FACET_WIDTH;
- Column 9: The number of elements of the Hilbert Basis, N_HILBERT_BASIS.

## 3. Combinatorial Properties

**Property 3.1.** *All the lattice points in the polytope defined by Definition 1.1 are boundary lattice points.*

*Proof.* Assuming there is an internal lattice point and the integer associated with the lattice point is M. Since all coordinates of an internal lattice point should be positive integers, or $x_{i,M} \geq 1$, minimum value of M would be achieved by setting $x_{i,M} = 1$. From Definition 1.1, we have

$$M = \prod_{i=1}^{n} p_i^{x_{i,M}} \geq \prod_{i=1}^{n} p_i > N$$

as any integer $N$ is smaller than the product of the all the prime numbers less than or equal to the number $N$ itself. The above equation implies that $N > M$, which contradicts $M \leq N$ as defined in Definition 1.1.

Therefore, the polytope defined by Definition 1.1 has no internal lattice points, and all lattice points are located on the boundary of the polytope. This type of polytope is called hollow polytope [7].

Property 3.1 is confirmed by the coefficients of the $h^*$ polynomials and the degrees of the polytopes representing natural numbers up to 66 given in Table 3.1, in which the degree of the polytope (or the degree of the $h^*$ vector), deg (P), is less than the dimension of the polytope, that is

$$\deg(P) < n = \dim(P)$$

The number of interior lattice points is related to the $n$-th coefficient of the $h^*$ polynomial by (Proposition 1.3.3, Gabriele Balletti [8])

$$h_n^* = |P^o \cap \mathbf{Z}^n|$$

where $P^o$ denotes the interior of the polytope $P$. We see, from the $h^*$ coefficients given in Table 3.1, that

$$h_n^* = 0$$

which shows that there is no interior point in the polytopes defined in Definition 1.1.

**Table 3.1. Coefficients of the h\* polynomials for polytopes representing natural numbers up to 66**

| N | dim(P) | Vol(P) | deg(P) | h₀* | h₁* | h₂* | h₃* | h₄* |
|---|---|---|---|---|---|---|---|---|
| 1 | 0 | 1 | | | | | | |
| 2 | 1 | 1 | 0 | 1 | | | | |
| 3 | 2 | 1 | 0 | 1 | | | | |
| 4 | 2 | 2 | 1 | 1 | 1 | | | |
| 5 | 3 | 2 | 1 | 1 | 1 | | | |
| 6 | 3 | 3 | 1 | 1 | 2 | | | |
| 7 | 4 | 3 | 1 | 1 | 2 | | | |
| 8 | 4 | 4 | 1 | 1 | 3 | | | |
| 9 | 4 | 6 | 2 | 1 | 4 | 1 | | |
| 10 | 4 | 8 | 2 | 1 | 5 | 2 | | |
| 11 | 5 | 8 | 2 | 1 | 5 | 2 | | |
| 12 | 5 | 9 | 2 | 1 | 6 | 2 | | |
| 13 | 6 | 9 | 2 | 1 | 6 | 2 | | |
| 14 | 6 | 11 | 2 | 1 | 7 | 3 | | |
| 15 | 6 | 14 | 2 | 1 | 8 | 5 | | |
| 16 | 6 | 16 | 2 | 1 | 9 | 6 | | |
| 17 | 7 | 16 | 2 | 1 | 9 | 6 | | |
| 18 | 7 | 18 | 2 | 1 | 10 | 7 | | |
| 19 | 8 | 18 | 2 | 1 | 10 | 7 | | |
| 20 | 8 | 20 | 2 | 1 | 11 | 8 | | |
| 21 | 8 | 24 | 2 | 1 | 12 | 11 | | |
| 22 | 8 | 28 | 2 | 1 | 13 | 14 | | |
| 23 | 9 | 28 | 2 | 1 | 13 | 14 | | |
| 24 | 9 | 29 | 2 | 1 | 14 | 14 | | |
| 25 | 9 | 38 | 3 | 1 | 15 | 20 | 2 | |
| 26 | 9 | 44 | 3 | 1 | 16 | 24 | 3 | |
| 27 | 9 | 52 | 3 | 1 | 17 | 29 | 5 | |
| 28 | 9 | 58 | 3 | 1 | 18 | 33 | 6 | |
| 29 | 10 | 58 | 3 | 1 | 18 | 33 | 6 | |
| 30 | 10 | 59 | 3 | 1 | 19 | 33 | 6 | |
| 31 | 11 | 59 | 3 | 1 | 19 | 33 | 6 | |
| 32 | 11 | 63 | 3 | 1 | 20 | 36 | 6 | |
| 33 | 11 | 73 | 3 | 1 | 21 | 42 | 9 | |
| 34 | 11 | 83 | 3 | 1 | 22 | 48 | 12 | |
| 35 | 11 | 100 | 3 | 1 | 23 | 56 | 20 | |
| 36 | 11 | 102 | 3 | 1 | 24 | 57 | 20 | |
| 37 | 12 | 102 | 3 | 1 | 24 | 57 | 20 | |
| 38 | 12 | 115 | 3 | 1 | 25 | 64 | 25 | |
| 39 | 12 | 131 | 3 | 1 | 26 | 72 | 32 | |
| 40 | 12 | 135 | 3 | 1 | 27 | 75 | 32 | |
| 41 | 13 | 135 | 3 | 1 | 27 | 75 | 32 | |
| 42 | 13 | 139 | 3 | 1 | 28 | 78 | 32 | |
| 43 | 14 | 139 | 3 | 1 | 28 | 78 | 32 | |
| 44 | 14 | 148 | 3 | 1 | 29 | 83 | 35 | |
| 45 | 14 | 159 | 3 | 1 | 30 | 89 | 39 | |
| 46 | 14 | 176 | 3 | 1 | 31 | 98 | 46 | |
| 47 | 15 | 176 | 3 | 1 | 31 | 98 | 46 | |
| 48 | 15 | 178 | 3 | 1 | 32 | 99 | 46 | |
| 49 | 15 | 216 | 3 | 1 | 33 | 114 | 68 | |
| 50 | 15 | 226 | 3 | 1 | 34 | 120 | 71 | |
| 51 | 15 | 252 | 3 | 1 | 35 | 131 | 85 | |
| 52 | 15 | 266 | 3 | 1 | 36 | 138 | 91 | |
| 53 | 16 | 266 | 3 | 1 | 36 | 138 | 91 | |
| 54 | 16 | 276 | 3 | 1 | 37 | 144 | 94 | |
| 55 | 16 | 334 | 4 | 1 | 38 | 158 | 123 | 14 |
| 56 | 16 | 343 | 4 | 1 | 39 | 162 | 127 | 14 |
| 57 | 16 | 380 | 4 | 1 | 40 | 175 | 147 | 17 |
| 58 | 16 | 420 | 4 | 1 | 41 | 188 | 169 | 21 |
| 59 | 17 | 420 | 4 | 1 | 41 | 188 | 169 | 21 |
| 60 | 17 | 422 | 4 | 1 | 42 | 189 | 169 | 21 |
| 61 | 18 | 422 | 4 | 1 | 42 | 189 | 169 | 21 |
| 62 | 18 | 462 | 4 | 1 | 43 | 202 | 191 | 25 |
| 63 | 18 | 478 | 4 | 1 | 44 | 211 | 197 | 25 |
| 64 | 18 | 492 | 4 | 1 | 45 | 218 | 203 | 25 |
| 65 | 18 | 563 | 4 | 1 | 46 | 235 | 241 | 40 |
| 66 | 18 | 578 | 4 | 1 | 47 | 242 | 247 | 41 |

The results given in Table 3.1 also confirm that the sum of the h* coefficients is equal to the relative volume of the polytope [8].

$$\sum_{i=0}^{n} h_i^* = Vol\ (P).$$

**Property 3.2.** *The coefficient $h_1^*$ of the polytope defined by Definition 1.1 for the natural number N is related to the dimension of the polytope by $h_1^* = N - n - 1$.*

*Proof.* According to Matthias Beck and Sinai Robins [9], if $P \in \mathbf{R}^n$ is a n-dimensional polytope and

$$L_P(t) := |tP \cap \mathbf{Z}^n| = \sum_{k=0}^{n} c_k t^k = \sum_{j=0}^{n} h_j^* \binom{t+n-j}{n}$$

Then

$$h_1^* = L_P(1) - n - 1.$$

Since the number of lattice points in the polytopes defined in Definition 1.1 for a natural number $N$ is equal to the number $N$ itself, we have

$$h_1^* = N - n - 1,$$

a relationship confirmed by the results shown in Table 3.1.

**Remark 3.1.** *The ratio of the number of vertices to the number of lattice points approaches 0.5 as N increases*

The ratio of the number of vertices to the number of lattice points is shown in Figure 3.1. Starting from 1 (for polytopes representing the numbers 1, 2, and 3), the ratio zigzags downward and becomes more stabilized around 0.5 as the number of vertices increases, although the ratio of a few data points goes below 0.5.

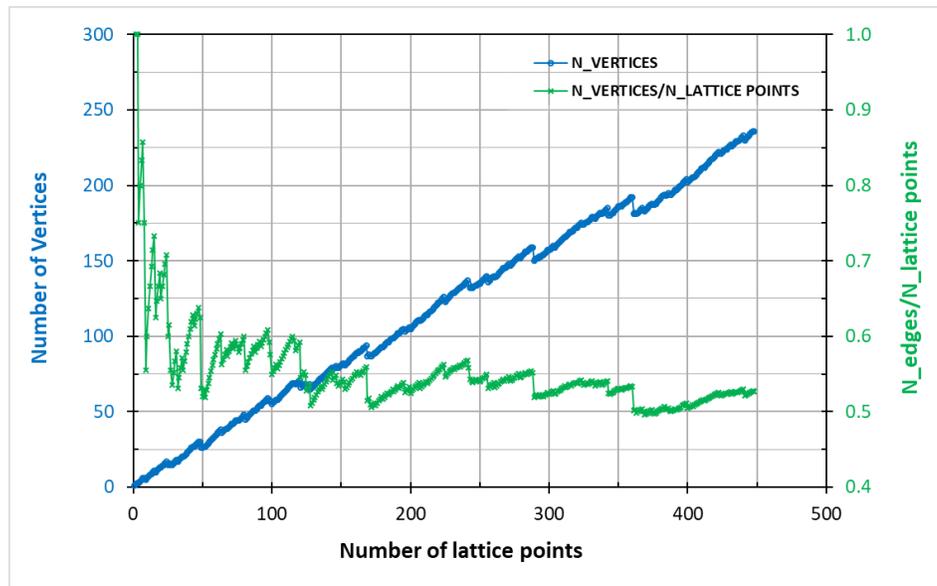

Figure 3.1 Number of vertices and N_Vertices/N_Lattice points ratio

**Remark 3.2.** *The F-vectors of the polytopes defined in Definition 1.1 satisfy the Euler–Poincaré–Schläfli formula:*

$$\sum_{i=-1}^{n}(-1)^i f_i = 0$$

where $f_{-1} = f_n = 1$.

The entries of the F-vectors of the polytopes for natural numbers up to 69 are given in Appendix B.

## 4. Geometric Properties

**Property 4.1.** *The Diameter of the polytope defined by Definition 1.1 is less than or equal to 2, or*

$$\delta(P) \leq 2.$$

The distance between two vertices of a polytope is the minimum number of edges in a path joining the two vertices. The diameter of a polytope is the greatest distance between two vertices of the polytope [10]. It can be seen from Appendix A that the diameter of all polytopes is 2, except that of the polytopes for $N = 2, 3, 4, 5, 9$, which are simplexes with diameter of 1. This means, any two vertices on the polytope representing a natural number are connected through at most two edges.

The diameters of the polytopes for the natural numbers up to 11 satisfy

$$\delta(P) = N_{facet} - n,$$

where $N_{facet}$ is the number of facets and $n$ is the dimension of the polytope. The values of the diameters of the polytopes given in Appendix A do not violate the Hirsch conjecture,

$$\delta(P) \leq N_{facet} - n,$$

although the Hirsch conjecture has been shown not true in all cases [11].

Here are some ideas regarding the proof of Property 4.1.

**Lemma 4.1.** *The lattice points associated with the $g_n$ consecutive integers, $p_n, \ldots, p_n + (g_n - 1)$ in the polytope representing the number $p_n + (g_n - 1)$ as defined in Definition 1.1 are all vertices of the polytope.*

**Lemma 4.2.** *The $g_n$ vectors from the origin to each of the $g_n$ vertices defined in Lemma 4.1 are linearly independent, and the $g_n$ vertices form a simplex of dimension $g_n - 1$.*

**Corollary 4.1.** *It follows from Lemma 4.1 and Lemma 4.2 that the $g_{n-1}$ vectors from the origin to each of the $g_{n-1}$ vertices associated with the $g_{n-1}$ consecutive integers, $p_{n-1}, \ldots, p_{n-1} + (g_{n-1} - 1)$ in the polytope representing the number $p_{n-1} + (g_{n-1} - 1)$ are linearly independent, and the $g_{n-1}$ vertices form a simplex of dimension $g_{n-1} - 1$.*

**Lemma 4.3**. *Combining the vertices in the $(g_{n-1} - 1)$-dimensional simplex defined in Corollary 4.1 with the vertex associated with the number $p_n$ forms a $(g_{n-1})$-dimensional pyramid.*

*Proof of Property 4.1.* Property 4.1 is true for $N = 1$. For a prime number $p_n$, any vertex in the pyramid Defined in Lemma 4.3 can be reached along the edge leading to the apex of the pyramid and the edge from the apex to the other vertex. Any vertex in the simplex defined in Lemma 4.2 is connected by another vertex in the simplex. If we constructed a new polytope by combining the vertexes on the pyramid and those on the simplex, any vertex can be reached by another vertex in the polytope through the two edges: a) from one of the two vertexes to the pyramid apex associated with the prime number $p_n$; and b) from the apex to the other vertex. Therefore, the diameter of the new polytope is less than or equal to 2.

**Remark 4.1.** *The volume of the polytope defined by Definition 1.1 decreases with N and asymptotically approaches zero as N approaches infinity.*

The volume given in Appendix A (Column 4) is the normalized the volume of the polytopes. The normalized volume of a lattice polytope has the advantage of being always an integer. The volume of a polytope, denoted as vol(P), is the normalized volume Vol(P) divided by the factorial of the dimension of the polytope.

$$vol(P) = \frac{Vol(P)}{n!}$$

where $n := \pi(N) = \dim(P)$, the dimension of polytope P.

The volume of the polytope representing natural numbers up to 448 is shown in Figure 4.1. In general, the volume decreases exponentially as the natural number increases.

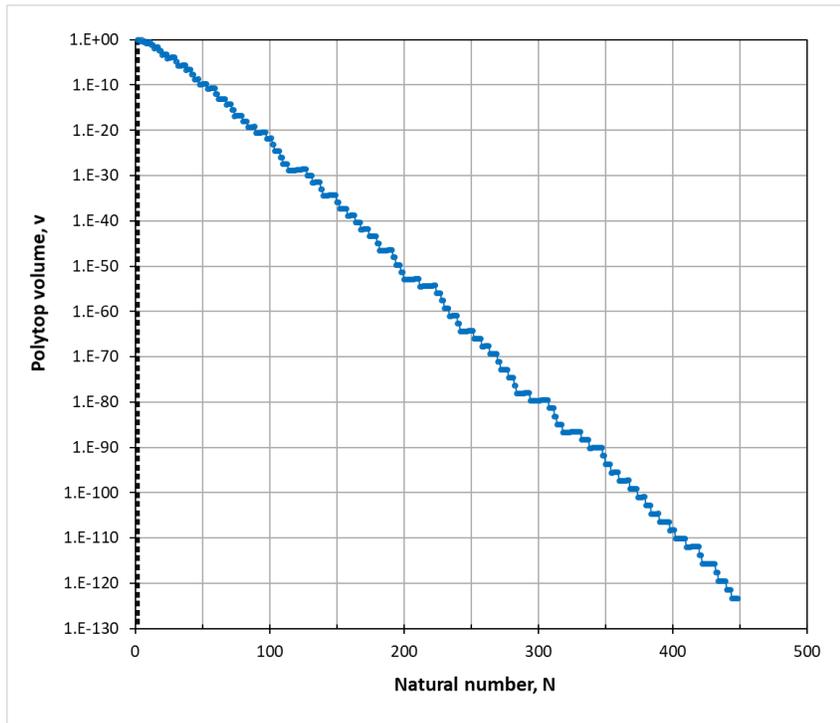

Figure 4.1 Volume of polytopes representing natural numbers up to 448

**Property 4.2.** *The volume of the polytope representing natural numbers is less than or equal to 1, or*

$$vol(P) \leq 1.$$

We know from Property 3.1 that the polytope defined by Definition 1.1 is hollow and does not have any interior lattice points. Let's determine the volume of this type of polytope by assuming that the n-dimensional polytope representing the natural number N has just one interior lattice point, L and $1 < M \leq N$.

We need the point L stay as close as possible to the origin so that the volume of the polytope is minimized. In doing so, the coordinates of L have to be 1 in all dimensions.

The (n-1)-dimensional affine plane in which the lattice point L is located satisfy the following equation

$$\sum_{i=1}^{n} \frac{x_i}{a_i} = 1$$

where $x_i$ is the Cartesian axis in *i-th* dimension and integer $a_i$ is the intercept on $x_i$ axis.

The volume of the n-dimensional simplex formed by the (n-1)-dimensional affine plane and the half spaces $x_i \geq 0$ is

$$V(\Delta_n) = \frac{\prod_{i=1}^{n} a_i}{n!}$$

Since the point L(1,1, ..., 1) is on the affine plane, the following condition should be satisfied.

$$\sum_{i=1}^{n} \frac{1}{a_i} = 1$$

It can be shown that the minimal volume of $n^n/n!$ is achieved at $a_i = n$.

In order to make point L an interior point of the n-dimensional simplex, the affine plane should be moved above the point L, which can be achieved by increasing one of the intercepts by 1 (from n to n + 1). In this case, the volume of the n-simplex increased to $(n + 1) n^{n-1}/n!$. Therefore, the minimal volume of the n-dimensional simplex representing the number, N, with one interior lattice point is

$$vol(P)_{min} = \frac{n^{n-1}(n+1)}{n!}$$

As $(n + 1) n^{n-1}/n!$ is always greater than 1 for $n \geq 1$, which is contradictory to the fact that all the polytope volumes listed in Appendix A are less than or equal to 1. Thus, none of the polytopes should have any internal lattice points, as far as the polytopes listed in Appendix A are concerned.

**Property 4.3.** *The volume of a polytope for a bigger natural number is larger than that for a smaller number having the same prime count as the bigger number.*

$$vol(P_N) > vol(P_M) \qquad if\ N > M\ and\ \pi(N) = \pi(M)$$

Increasing a natural number before reaching the next prime number is done geometrically in constructing the polytope by adding a new lattice point to the existing polytope of the same dimension. Essentially what Property 4.3 says is that such a lattice point will always be added somewhere outside (not inside or on the boundary of) the existing polytope and this lattice point will become a new vertex of the polytope. It is obvious that, if the dimension of the polytope remains the same, adding a vertex outside an existing polytope creates a new polytope with bigger Euclidian volume, $Vol(P)$, resulting in a bigger polytope volume, defined as $vol(P) = Vol(P)/n!$.

**Property 4.4.** *The polytope volume representing a prime number, $p$, is $\pi(p)$ times smaller than that representing the previous number, $p - 1$.*

$$vol(P_p) = \frac{vol(P_{p-1})}{\pi(p)}$$

*Proof.* Let's consider a number, $p - 1$. The volume of the polytope $P_{p-1} \in \mathbf{R}^{n-1}$ that represents the number $p - 1$ is $vol(P_{p-1})$. If the next number p is a prime number, the polytope for the number $p$ is constructed by using the polytope $P_{p-1}$ as its base and adding a new dimension perpendicular to the base polytope $P_{p-1}$. From Definition 1.1, the number $p$ is expressed by

$$p = \prod_{i=1}^{n} p_i^{x_{i,p}}$$

Since $p$ is a prime number and $p = p_n$, the exponents $x_{i,p} = 0$ except that $x_{n,p} = 1$. The coordinates of the newly added vertex by the prime number $p$ is p(0, 0, 0, ... , 1) $\in \mathbf{R}^n$. The unit vector from the origin to the vertex $p$ is perpendicular to the base polytope, and the volume of the polytope constructed, $P_p = conv\{P_{p-1} \cup \{e_n\}\}$, which is a pyramid, can be calculated by multiplying the volume of the base polytope by the height of this unit vector added and divided by the dimension number of the new dimension, $\pi(p)$.

**Property 4.5.** *The polytope volume for the number before a prime number decrease as prime count increases*

Another observation can be made from the Appendix A is that, as the prime count increase the polytope volume of the number before a prime number decreases, except for the cases from 1 to 2 and from 1 (or 2) to 4, in which the polytope volume stays the same.

If M is the largest number among all the consecutive numbers sharing the same prime count $\pi(M)$, and N is the largest numbers among those sharing the same prime count $\pi(N)$, then

$$vol(P_N) \leq vol(P_M) \qquad\qquad if\ M < N$$

where $M = p_s - 1$, and $N = p_t - 1$, in which s and t are two positive integers and $s < t$. $p_s$ and $p_t$ are the *s-th* and *t-th* prime numbers.

$v(N) = v(M)$ holds only for three cases: i) $M = 1, N = 2$; ii) $M = 1, N = 4$, and iii) $M = 2, N = 4$.

**Property 4.6.** *The polytope volume for a smaller prime number is larger that of a bigger prime number*

For any two prime numbers, $P_1$ and $P_2$ with $p_1 < p_2$, then $v(p_1) > v(p_2)$.

Based on the Properties 3.4.2 through 3.4.5, we have the following conjecture

**Conjecture 4.1.** *The volume of polytope defined by Definition 1.1 asymptotically approaches zero as the natural number N approaches infinity, or*

$$\lim_{N\to\infty} v(N) = 0.$$

It is known since Euclidean (c. 300 BC) that there are infinite number of prime numbers. Let us see what it means from the perspective of polytope volume if there were finite number of prime numbers.

Assuming there exist only *m* prime numbers, 2, 3, 5, ..., $p_m$. Since all the number greater than $p_m$ would share the same prime count $\pi(p_m)$, all the polytopes for numbers greater than $p_m$ would be *m*-dimensional. Thus, if the number of prime numbers were finite,

- Polytopes representing any number greater than or equal to $\prod_{i=1}^{m} p_i$ would have an interior lattice point (1, 1, 1, ..., 1) in $R^m$;

- The volume of the polytope representing any number greater than $p_m$ would increase monotonously as the number increases; and

- The sum of the volumes of the polytopes representing all the natural numbers approach infinity.

However, none of the above is observed from the results shown in Appendix A.

**Remark 4.2.** *The normalized volume of polytopes representing a natural number N, denoted as $Vol(N)$, given in Appendix A can be correlated in terms of the number of vertices, $N_v$, and the dimension of the polytope, n, by the following equation;*

$$Vol(N) = 2(lnN_v)^{\sqrt{n}}$$

The normalized volumes are compared with the correlated values calculated using the equation above as a function of N in Figure 4.2.

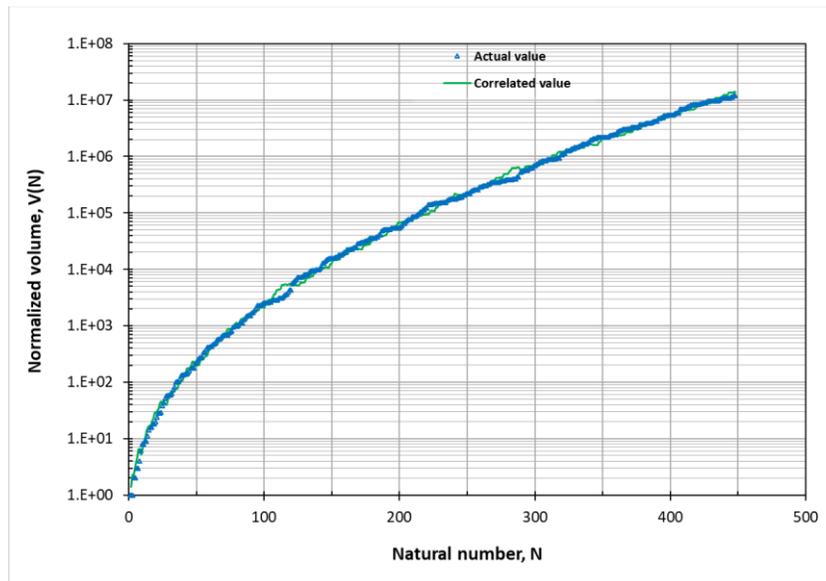

Figure 4.2 Actual and correlated values of normalized volumes

**Property 4.7.** *The sum of the volumes of the polytopes representing all the natural numbers is equal to approximately 5.42915412.*

Based on the normalized volumes listed in Appendix A, we calculated the sum of the volumes of the polytopes representing the number from 1 to N by the following equation.

$$\text{vol}_{sum}(P_N) = \sum_{M=1}^{N} \text{vol}(P_M) = \sum_{M=1}^{N} \frac{\text{Vol}(M)}{\pi(M)!}$$

where $Vol(M)$ is the normalized volume of the polytope representing the natural number M.

The values of $v_{sum}(N)$ with 150 digits of decimal points for N up to 448 were calculated. The sum of the polytope volume converges to a constant of approximately 5.42915412 as N increases. Let's name this constant as ρ,

$$\rho := \sum_{N=0}^{\infty} v(N)$$

The first 123 digits of the constant ρ is

5.429154123959951333182133057457967726075977448089132317290226047820875978258078787063403516637498516150454315367601796 42454.

The sum of the volumes representing all prime numbers is equal to e, the base of the natural logarithm (see Section 5). Thus, the ratio of the polytope volume of all the natural numbers to that of all the prime numbers is ρ/e. The value of ρ/e equals approximately to 1.997274185156, which is close to 2.

## 5. Polytopes Representing a Subset of Natural Numbers

Definition 1.1 can be modified to cover a subset of natural numbers. In this case, the polytope representing the same number could be different from that when all the natural numbers are considered. If we consider only prime numbers, the polytope representing prime numbers will become simplexes. The polytope representing the number 7, for example, is a 3-simplex with a volume of 1/24. However, if we consider all the natural numbers, the polytope representing the number 7 is a 4-polytope with a volume of 1/8.

In both scenarios, the polytope is constructed based on the principle that each of the prime numbers represented as a dimension in the Euclidian space. With the same principle, all odd numbers could be represented by polytopes. However, if prime numbers are used as the basis of the polytopes, the number "2" must be included in the subset simply because 2 is a prime number. Likewise, all the even numbers and prime numbers as a group can also be represented with prime number-based polytopes.

Table 5.1 gives the volumes of polytopes representing the first few numbers for different subsets of natural numbers. The sum of the polytope volumes is given in the last row in Table 5.1. The sum of the volumes of a subset of natural numbers is unique and can be considered as the "signature" of the subset of the natural numbers.

The sum of the volumes for "1, even numbers and odd prime numbers" is

$\rho_1$ =4.9596791660088175666106133110947935933119252633787730402179565130724230464445 02259...

and the ratio $\rho_1/e$=1.824563999980...

The sum of the volumes for "2 and odd numbers" is
$\rho_2$ =2.8160992749048251156406996681472257475482767405230146224656363186575797069 62 24029...

and the ratio $\rho_2/e$=1.035985027535...

$\{N^2\}, \{N^3\}$, ect. are examples of other subsets of natural numbers. The polytopes representing these types of natural numbers are also given in Table 5.1.

Table 5.1 Volumes (and dimensions) of polytopes representing numbers for different subsets of natural numbers

| Subset | 1 and prime numbers | Natural number | 1, even numbers and odd primes | 2 and odd numbers | $N^2$ | $N^3$ |
|---|---|---|---|---|---|---|
| $vol(P_1)$ | 1 (0)* | | | | | |
| $vol(P_2)$ | 1 (1) | | | | | |
| $vol(P_3)$ | 1/2 (2) | | | | | |
| $vol(P_4)$ | | 1 (2) | | | 2 (1) | |
| $vol(P_5)$ | 1/6 (3) | 1/3 (3) | | 1/6 (3) | | |
| $vol(P_6)$ | | 1/2 (3) | | | | |
| $vol(P_7)$ | 1/24 (4) | 1/8 (4) | | 1/24 (4) | | |
| $vol(P_8)$ | | 1/6 (4) | | | | 3(1) |
| $vol(P_9)$ | | 1/4 (4) | | 1/12 (4) | 2 (2) | |
| $vol(P_{10})$ | | 1/3 (4) | 5/24 (4) | | | |
| $vol(P_{11})$ | 1/120 (5) | 1/15 (5) | 1/24 (5) | 1/60 (5) | | |
| $\sum_{M=1}^{\infty} vol(P_M)$ | e=2.71828... | 5.4291541239... | 4.9596791660... | 2.8160992749... | 40.492978173... | 243.5229862... |
| $\sum_{M=1}^{\infty} vol(P_M)/e$ | 1 | 1.9972741851... | 1.82456399998... | 1.03598502753... | | |

* the number given in the parenthesis is the dimension of the polytope.

In general, any subset of the natural numbers, *S ϵ* **N,** can be used as the basis in constructing the polytopes representing another subset of natural numbers, *T ϵ* **N,** as long as every number in T can be expressed uniquely as the product of the elements in S, just like the canonical representation[12] of natural numbers by prime numbers.

**Conclusions**

A representation method was proposed to represent natural numbers by lattice polytopes, and the characteristics of the polytopes were studied using *Polymak* and *Qhull* software. It was observed from the results that these polytopes do not have any interior lattice points, and any two vertices are connected through no more than two edges. The volume of the polytopes approaches zero as the natural number approaches infinity and the sum of the volumes of the polytopes representing all the natural numbers converges to a constant.

Yaping Lu: luyaping1@yahoo.com

Shufang Deng: sfangd@163.com

# Appendix A: Properties of polytopes representing natural numbers up to 448

| N | Diam. δ | Dim. n | Normalized volume, Vol | Number of vertices | Number of edges | Number of facets | Maximal facet width | Number of Hilbert basis |
|---|---------|--------|------------------------|--------------------|-----------------|------------------|---------------------|-------------------------|
| 1 |   | 0 | 1 | 1 |   |   |   |   |
| 2 | 1 | 1 | 1 | 2 | 1 | 2 | 1 | 2 |
| 3 | 1 | 2 | 1 | 3 | 3 | 3 | 1 | 3 |
| 4 | 1 | 2 | 2 | 3 | 3 | 3 | 2 | 4 |
| 5 | 1 | 3 | 2 | 4 | 6 | 4 | 2 | 5 |
| 6 | 2 | 3 | 3 | 5 | 8 | 5 | 2 | 6 |
| 7 | 2 | 4 | 3 | 6 | 13 | 6 | 2 | 7 |
| 8 | 2 | 4 | 4 | 6 | 13 | 6 | 3 | 8 |
| 9 | 1 | 4 | 6 | 5 | 10 | 5 | 6 | 9 |
| 10 | 2 | 4 | 8 | 6 | 13 | 6 | 6 | 10 |
| 11 | 2 | 5 | 8 | 7 | 19 | 7 | 6 | 11 |
| 12 | 2 | 5 | 9 | 8 | 23 | 8 | 4 | 12 |
| 13 | 2 | 6 | 9 | 9 | 31 | 9 | 4 | 13 |
| 14 | 2 | 6 | 11 | 10 | 35 | 9 | 4 | 14 |
| 15 | 2 | 6 | 14 | 11 | 40 | 10 | 4 | 15 |
| 16 | 2 | 6 | 16 | 10 | 36 | 10 | 4 | 16 |
| 17 | 2 | 7 | 16 | 11 | 46 | 11 | 4 | 17 |
| 18 | 2 | 7 | 18 | 12 | 51 | 11 | 8 | 18 |
| 19 | 2 | 8 | 18 | 13 | 63 | 12 | 8 | 19 |
| 20 | 2 | 8 | 20 | 13 | 63 | 12 | 8 | 20 |
| 21 | 2 | 8 | 24 | 14 | 70 | 13 | 8 | 21 |
| 22 | 2 | 8 | 28 | 15 | 76 | 13 | 8 | 22 |
| 23 | 2 | 9 | 28 | 16 | 91 | 14 | 8 | 23 |
| 24 | 2 | 9 | 29 | 17 | 97 | 13 | 5 | 24 |
| 25 | 2 | 9 | 38 | 15 | 83 | 13 | 10 | 25 |
| 26 | 2 | 9 | 44 | 16 | 89 | 13 | 10 | 26 |
| 27 | 2 | 9 | 52 | 15 | 83 | 13 | 18 | 27 |
| 28 | 2 | 9 | 58 | 15 | 84 | 14 | 18 | 29 |
| 29 | 2 | 10 | 58 | 16 | 99 | 15 | 18 | 30 |
| 30 | 2 | 10 | 59 | 17 | 108 | 15 | 9 | 30 |
| 31 | 2 | 11 | 59 | 18 | 125 | 16 | 9 | 31 |
| 32 | 2 | 11 | 63 | 17 | 116 | 16 | 15 | 32 |
| 33 | 2 | 11 | 73 | 18 | 126 | 17 | 15 | 33 |
| 34 | 2 | 11 | 83 | 19 | 135 | 17 | 15 | 34 |
| 35 | 2 | 11 | 100 | 20 | 149 | 22 | 15 | 35 |
| 36 | 2 | 11 | 102 | 20 | 148 | 21 | 10 | 36 |
| 37 | 2 | 12 | 102 | 21 | 168 | 22 | 10 | 37 |
| 38 | 2 | 12 | 115 | 22 | 178 | 22 | 10 | 38 |
| 39 | 2 | 12 | 131 | 23 | 189 | 22 | 10 | 39 |
| 40 | 2 | 12 | 135 | 24 | 200 | 21 | 10 | 40 |
| 41 | 2 | 13 | 135 | 25 | 224 | 22 | 10 | 41 |
| 42 | 2 | 13 | 139 | 26 | 233 | 21 | 10 | 42 |
| 43 | 2 | 14 | 139 | 27 | 259 | 22 | 10 | 43 |
| 44 | 2 | 14 | 148 | 27 | 259 | 22 | 10 | 44 |
| 45 | 2 | 14 | 159 | 28 | 275 | 26 | 12 | 45 |
| 46 | 2 | 14 | 176 | 29 | 287 | 26 | 12 | 46 |
| 47 | 2 | 15 | 176 | 30 | 316 | 27 | 12 | 47 |
| 48 | 2 | 15 | 178 | 30 | 314 | 26 | 12 | 48 |
| 49 | 2 | 15 | 216 | 26 | 259 | 24 | 22 | 49 |
| 50 | 2 | 15 | 226 | 26 | 255 | 22 | 34 | 50 |
| 51 | 2 | 15 | 252 | 27 | 267 | 22 | 34 | 51 |
| 52 | 2 | 15 | 266 | 27 | 268 | 22 | 34 | 52 |
| 53 | 2 | 16 | 266 | 28 | 295 | 23 | 34 | 53 |
| 54 | 2 | 16 | 276 | 29 | 308 | 22 | 22 | 55 |
| 55 | 2 | 16 | 334 | 30 | 329 | 25 | 22 | 56 |
| 56 | 2 | 16 | 343 | 31 | 346 | 29 | 36 | 56 |
| 57 | 2 | 16 | 380 | 32 | 360 | 29 | 36 | 57 |
| 58 | 2 | 16 | 420 | 33 | 374 | 29 | 36 | 58 |
| 59 | 2 | 17 | 420 | 34 | 407 | 30 | 36 | 59 |
| 60 | 2 | 17 | 422 | 35 | 422 | 30 | 18 | 60 |
| 61 | 2 | 18 | 422 | 36 | 457 | 31 | 18 | 61 |
| 62 | 2 | 18 | 462 | 37 | 472 | 31 | 18 | 62 |
| 63 | 2 | 18 | 478 | 38 | 494 | 33 | 16 | 63 |
| 64 | 2 | 18 | 492 | 36 | 461 | 33 | 18 | 64 |
| 65 | 2 | 18 | 563 | 37 | 483 | 32 | 18 | 65 |
| 66 | 2 | 18 | 578 | 38 | 495 | 32 | 18 | 66 |
| 67 | 2 | 19 | 578 | 39 | 533 | 33 | 18 | 67 |
| 68 | 2 | 19 | 610 | 39 | 535 | 35 | 18 | 68 |
| 69 | 2 | 19 | 671 | 40 | 552 | 35 | 18 | 69 |
| 70 | 2 | 19 | 684 | 41 | 573 | 34 | 18 | 70 |
| 71 | 2 | 20 | 684 | 42 | 614 | 35 | 18 | 71 |
| 72 | 2 | 20 | 687 | 42 | 612 | 33 | 12 | 72 |
| 73 | 2 | 21 | 687 | 43 | 654 | 34 | 12 | 73 |
| 74 | 2 | 21 | 746 | 44 | 673 | 34 | 12 | 74 |
| 75 | 2 | 21 | 773 | 44 | 671 | 36 | 24 | 75 |
| 76 | 2 | 21 | 808 | 44 | 671 | 36 | 24 | 76 |
| 77 | 2 | 21 | 936 | 45 | 700 | 41 | 24 | 77 |
| 78 | 2 | 21 | 963 | 46 | 716 | 39 | 24 | 78 |
| 79 | 2 | 22 | 963 | 47 | 762 | 40 | 24 | 79 |
| 80 | 2 | 22 | 972 | 48 | 783 | 40 | 19 | 80 |
| 81 | 2 | 22 | 1022 | 45 | 723 | 40 | 32 | 83 |
| 82 | 2 | 22 | 1103 | 46 | 742 | 40 | 32 | 84 |
| 83 | 2 | 23 | 1103 | 47 | 788 | 41 | 32 | 85 |
| 84 | 2 | 23 | 1110 | 48 | 811 | 40 | 32 | 86 |
| 85 | 2 | 23 | 1265 | 49 | 839 | 42 | 32 | 86 |
| 86 | 2 | 23 | 1361 | 50 | 859 | 42 | 32 | 87 |
| 87 | 2 | 23 | 1481 | 51 | 883 | 42 | 32 | 88 |
| 88 | 2 | 23 | 1509 | 51 | 884 | 41 | 40 | 90 |
| 89 | 2 | 24 | 1509 | 52 | 935 | 42 | 40 | 91 |
| 90 | 2 | 24 | 1511 | 53 | 957 | 41 | 20 | 90 |
| 91 | 2 | 24 | 1690 | 54 | 991 | 43 | 20 | 91 |
| 92 | 2 | 24 | 1753 | 54 | 989 | 43 | 20 | 92 |
| 93 | 2 | 24 | 1896 | 55 | 1014 | 43 | 20 | 93 |
| 94 | 2 | 24 | 2035 | 56 | 1037 | 43 | 20 | 94 |
| 95 | 2 | 24 | 2292 | 57 | 1067 | 43 | 20 | 95 |
| 96 | 2 | 24 | 2306 | 58 | 1086 | 44 | 20 | 96 |
| 97 | 2 | 25 | 2306 | 59 | 1144 | 45 | 20 | 97 |
| 98 | 2 | 25 | 2374 | 58 | 1116 | 49 | 40 | 98 |
| 99 | 2 | 25 | 2477 | 57 | 1101 | 55 | 40 | 99 |
| 100 | 2 | 25 | 2516 | 55 | 1054 | 57 | 40 | 101 |
| 101 | 2 | 26 | 2516 | 56 | 1109 | 58 | 40 | 102 |
| 102 | 2 | 26 | 2565 | 57 | 1128 | 59 | 40 | 103 |
| 103 | 2 | 27 | 2565 | 58 | 1185 | 60 | 40 | 104 |
| 104 | 2 | 27 | 2616 | 58 | 1186 | 56 | 40 | 105 |
| 105 | 2 | 27 | 2644 | 59 | 1219 | 55 | 40 | 105 |
| 106 | 2 | 27 | 2819 | 60 | 1245 | 55 | 40 | 106 |
| 107 | 2 | 28 | 2819 | 61 | 1305 | 56 | 40 | 107 |
| 108 | 2 | 28 | 2827 | 62 | 1329 | 53 | 26 | 108 |
| 109 | 2 | 29 | 2827 | 63 | 1391 | 54 | 26 | 109 |
| 110 | 2 | 29 | 2866 | 64 | 1416 | 52 | 26 | 110 |
| 111 | 2 | 29 | 3073 | 65 | 1445 | 52 | 26 | 111 |
| 112 | 2 | 29 | 3104 | 66 | 1480 | 53 | 34 | 112 |
| 113 | 2 | 30 | 3104 | 67 | 1546 | 54 | 34 | 113 |
| 114 | 2 | 30 | 3139 | 68 | 1568 | 50 | 34 | 114 |
| 115 | 2 | 30 | 3454 | 69 | 1606 | 54 | 34 | 115 |
| 116 | 2 | 30 | 3588 | 69 | 1607 | 54 | 34 | 116 |
| 117 | 2 | 30 | 3683 | 68 | 1587 | 54 | 34 | 117 |
| 118 | 2 | 30 | 3909 | 69 | 1616 | 54 | 34 | 118 |
| 119 | 2 | 30 | 4321 | 70 | 1655 | 57 | 34 | 119 |
| 120 | 2 | 30 | 4326 | 71 | 1680 | 55 | 21 | 120 |
| 121 | 2 | 30 | 5604 | 66 | 1513 | 55 | 42 | 122 |
| 122 | 2 | 30 | 5954 | 67 | 1541 | 55 | 42 | 123 |
| 123 | 2 | 30 | 6386 | 68 | 1571 | 55 | 42 | 124 |
| 124 | 2 | 30 | 6612 | 68 | 1572 | 55 | 42 | 125 |
| 125 | 2 | 30 | 7268 | 66 | 1523 | 68 | 54 | 129 |
| 126 | 2 | 30 | 7300 | 67 | 1550 | 63 | 54 | 131 |
| 127 | 2 | 31 | 7300 | 68 | 1617 | 64 | 54 | 132 |
| 128 | 2 | 31 | 7422 | 65 | 1529 | 66 | 42 | 135 |
| 129 | 2 | 31 | 7940 | 66 | 1559 | 66 | 42 | 136 |
| 130 | 2 | 31 | 8032 | 67 | 1587 | 61 | 42 | 135 |
| 131 | 2 | 32 | 8032 | 68 | 1654 | 62 | 42 | 136 |
| 132 | 2 | 32 | 8044 | 69 | 1688 | 65 | 42 | 133 |
| 133 | 2 | 32 | 8816 | 70 | 1730 | 64 | 42 | 134 |
| 134 | 2 | 32 | 9362 | 71 | 1761 | 64 | 42 | 135 |
| 135 | 2 | 32 | 9460 | 72 | 1799 | 65 | 63 | 138 |
| 136 | 2 | 32 | 9602 | 72 | 1800 | 66 | 63 | 140 |
| 137 | 2 | 33 | 9602 | 73 | 1872 | 67 | 63 | 141 |
| 138 | 2 | 33 | 9740 | 74 | 1897 | 67 | 63 | 142 |
| 139 | 2 | 34 | 9740 | 75 | 1971 | 68 | 63 | 143 |
| 140 | 2 | 34 | 9783 | 76 | 2009 | 71 | 90 | 144 |
| 141 | 2 | 34 | 10393 | 77 | 2043 | 71 | 90 | 145 |
| 142 | 2 | 34 | 10997 | 78 | 2076 | 71 | 90 | 146 |
| 143 | 2 | 34 | 12732 | 79 | 2142 | 92 | 90 | 146 |
| 144 | 2 | 34 | 12744 | 78 | 2107 | 89 | 90 | 147 |
| 145 | 2 | 34 | 14014 | 79 | 2153 | 92 | 90 | 147 |
| 146 | 2 | 34 | 14822 | 80 | 2187 | 92 | 90 | 148 |
| 147 | 2 | 34 | 15175 | 79 | 2145 | 98 | 54 | 152 |
| 148 | 2 | 34 | 15671 | 79 | 2147 | 98 | 54 | 155 |
| 149 | 2 | 35 | 15671 | 80 | 2226 | 99 | 54 | 156 |
| 150 | 2 | 35 | 15692 | 81 | 2261 | 94 | 62 | 154 |
| 151 | 2 | 36 | 15692 | 82 | 2342 | 95 | 62 | 155 |
| 152 | 2 | 36 | 15825 | 82 | 2341 | 83 | 62 | 156 |
| 153 | 2 | 36 | 16479 | 81 | 2316 | 95 | 62 | 158 |
| 154 | 2 | 36 | 16588 | 82 | 2361 | 102 | 48 | 156 |
| 155 | 2 | 36 | 17932 | 83 | 2407 | 102 | 48 | 157 |
| 156 | 2 | 36 | 17986 | 84 | 2439 | 98 | 48 | 158 |
| 157 | 2 | 37 | 17986 | 85 | 2523 | 99 | 48 | 159 |
| 158 | 2 | 37 | 18967 | 86 | 2561 | 99 | 48 | 160 |
| 159 | 2 | 37 | 20117 | 87 | 2600 | 99 | 48 | 161 |
| 160 | 2 | 37 | 20154 | 88 | 2634 | 94 | 48 | 162 |
| 161 | 2 | 37 | 22086 | 89 | 2687 | 98 | 48 | 163 |
| 162 | 2 | 37 | 22301 | 89 | 2680 | 106 | 74 | 164 |
| 163 | 2 | 38 | 22301 | 90 | 2769 | 107 | 74 | 165 |
| 164 | 2 | 38 | 22986 | 90 | 2767 | 107 | 74 | 166 |
| 165 | 2 | 38 | 23100 | 91 | 2819 | 111 | 74 | 165 |
| 166 | 2 | 38 | 24297 | 92 | 2859 | 111 | 74 | 166 |
| 167 | 2 | 39 | 24297 | 93 | 2951 | 112 | 74 | 167 |
| 168 | 2 | 39 | 24308 | 94 | 2990 | 110 | 46 | 168 |
| 169 | 2 | 39 | 28458 | 87 | 2660 | 93 | 74 | 176 |
| 170 | 2 | 39 | 28784 | 88 | 2694 | 93 | 74 | 177 |
| 171 | 2 | 39 | 29642 | 87 | 2668 | 93 | 74 | 179 |
| 172 | 2 | 39 | 30408 | 87 | 2666 | 93 | 74 | 181 |
| 173 | 2 | 40 | 30408 | 88 | 2753 | 94 | 74 | 182 |
| 174 | 2 | 40 | 30864 | 89 | 2788 | 95 | 74 | 183 |
| 175 | 2 | 40 | 31426 | 89 | 2803 | 87 | 102 | 188 |
| 176 | 2 | 40 | 31654 | 90 | 2849 | 93 | 102 | 189 |
| 177 | 2 | 40 | 33524 | 91 | 2892 | 93 | 102 | 190 |
| 178 | 2 | 40 | 35166 | 92 | 2932 | 93 | 102 | 191 |
| 179 | 2 | 41 | 35166 | 93 | 3024 | 94 | 102 | 192 |
| 180 | 2 | 41 | 35176 | 93 | 3024 | 92 | 58 | 192 |
| 181 | 2 | 42 | 35176 | 94 | 3117 | 93 | 58 | 193 |
| 182 | 2 | 42 | 35339 | 95 | 3173 | 100 | 46 | 189 |
| 183 | 2 | 42 | 37234 | 96 | 3217 | 100 | 46 | 190 |
| 184 | 2 | 42 | 37542 | 96 | 3214 | 92 | 46 | 191 |
| 185 | 2 | 42 | 40635 | 97 | 3269 | 92 | 46 | 191 |
| 186 | 2 | 42 | 41071 | 98 | 3306 | 92 | 46 | 192 |
| 187 | 2 | 42 | 48657 | 99 | 3383 | 114 | 46 | 193 |
| 188 | 2 | 42 | 50005 | 99 | 3382 | 114 | 46 | 195 |
| 189 | 2 | 42 | 50376 | 100 | 3432 | 107 | 48 | 197 |
| 190 | 2 | 42 | 50840 | 101 | 3470 | 106 | 48 | 198 |
| 191 | 2 | 43 | 50840 | 102 | 3571 | 107 | 48 | 199 |
| 192 | 2 | 43 | 50929 | 102 | 3565 | 113 | 48 | 200 |
| 193 | 2 | 44 | 50929 | 103 | 3667 | 114 | 48 | 201 |
| 194 | 2 | 44 | 53343 | 104 | 3712 | 114 | 48 | 202 |
| 195 | 2 | 44 | 53605 | 105 | 3776 | 122 | 68 | 198 |
| 196 | 2 | 44 | 54047 | 103 | 3670 | 115 | 48 | 200 |
| 197 | 2 | 45 | 54047 | 104 | 3773 | 116 | 48 | 201 |
| 198 | 2 | 45 | 54156 | 105 | 3819 | 111 | 60 | 202 |
| 199 | 2 | 46 | 54156 | 106 | 3924 | 112 | 60 | 203 |
| 200 | 2 | 46 | 54239 | 105 | 3878 | 106 | 76 | 205 |
| 201 | 2 | 46 | 57129 | 106 | 3925 | 106 | 76 | 206 |
| 202 | 2 | 46 | 59735 | 107 | 3973 | 106 | 76 | 207 |
| 203 | 2 | 46 | 65634 | 108 | 4039 | 118 | 76 | 208 |
| 204 | 2 | 46 | 65770 | 109 | 4078 | 108 | 76 | 209 |
| 205 | 2 | 46 | 71106 | 110 | 4137 | 108 | 76 | 209 |
| 206 | 2 | 46 | 74302 | 111 | 4186 | 108 | 76 | 210 |
| 207 | 2 | 46 | 76136 | 110 | 4152 | 108 | 76 | 213 |
| 208 | 2 | 46 | 76503 | 111 | 4215 | 115 | 78 | 214 |
| 209 | 2 | 46 | 85347 | 112 | 4296 | 127 | 78 | 215 |
| 210 | 2 | 46 | 85356 | 113 | 4340 | 122 | 39 | 210 |
| 211 | 2 | 47 | 85356 | 114 | 4453 | 123 | 39 | 211 |
| 212 | 2 | 47 | 87354 | 114 | 4455 | 123 | 39 | 212 |
| 213 | 2 | 47 | 91932 | 115 | 4505 | 123 | 39 | 213 |
| 214 | 2 | 47 | 95953 | 116 | 4557 | 123 | 39 | 214 |
| 215 | 2 | 47 | 102747 | 117 | 4617 | 123 | 39 | 215 |
| 216 | 2 | 47 | 102808 | 118 | 4658 | 118 | 39 | 216 |
| 217 | 2 | 47 | 112222 | 119 | 4727 | 118 | 39 | 217 |
| 218 | 2 | 47 | 117088 | 120 | 4780 | 118 | 39 | 218 |
| 219 | 2 | 47 | 123031 | 121 | 4831 | 118 | 39 | 219 |
| 220 | 2 | 47 | 123385 | 122 | 4885 | 113 | 39 | 220 |
| 221 | 2 | 47 | 140750 | 123 | 4982 | 133 | 39 | 221 |
| 222 | 2 | 47 | 142242 | 124 | 5023 | 134 | 39 | 222 |
| 223 | 2 | 48 | 142242 | 125 | 5147 | 135 | 39 | 223 |
| 224 | 2 | 48 | 142526 | 126 | 5202 | 132 | 39 | 224 |
| 225 | 2 | 48 | 143227 | 123 | 5033 | 131 | 30 | 225 |
| 226 | 2 | 48 | 149291 | 124 | 5085 | 131 | 30 | 226 |
| 227 | 2 | 49 | 149291 | 125 | 5209 | 132 | 30 | 227 |
| 228 | 2 | 49 | 149498 | 126 | 5243 | 123 | 30 | 228 |
| 229 | 2 | 50 | 149498 | 127 | 5369 | 124 | 30 | 229 |
| 230 | 2 | 50 | 150648 | 128 | 5410 | 120 | 30 | 230 |
| 231 | 2 | 50 | 151625 | 129 | 5479 | 122 | 54 | 231 |
| 232 | 2 | 50 | 153563 | 129 | 5481 | 125 | 54 | 232 |
| 233 | 2 | 51 | 153563 | 130 | 5610 | 126 | 54 | 233 |
| 234 | 2 | 51 | 153786 | 131 | 5668 | 120 | 54 | 234 |
| 235 | 2 | 51 | 163843 | 132 | 5738 | 120 | 54 | 235 |
| 236 | 2 | 51 | 167803 | 132 | 5738 | 120 | 54 | 236 |
| 237 | 2 | 51 | 176134 | 133 | 5798 | 120 | 54 | 237 |
| 238 | 2 | 51 | 177171 | 134 | 5842 | 123 | 54 | 238 |
| 239 | 2 | 52 | 177171 | 135 | 5976 | 124 | 54 | 239 |
| 240 | 2 | 52 | 177185 | 136 | 6019 | 122 | 44 | 240 |
| 241 | 2 | 53 | 177185 | 137 | 6155 | 123 | 44 | 241 |
| 242 | 2 | 53 | 179658 | 135 | 6020 | 133 | 54 | 244 |
| 243 | 2 | 53 | 183115 | 132 | 5861 | 151 | 90 | 249 |
| 244 | 2 | 53 | 187171 | 132 | 5860 | 151 | 90 | 251 |
| 245 | 2 | 53 | 189778 | 132 | 5854 | 157 | 140 | 253 |
| 246 | 2 | 53 | 191706 | 133 | 5902 | 157 | 140 | 254 |
| 247 | 2 | 53 | 203552 | 134 | 6000 | 154 | 140 | 255 |
| 248 | 2 | 53 | 205479 | 134 | 6001 | 149 | 140 | 257 |
| 249 | 2 | 53 | 214788 | 135 | 6063 | 149 | 140 | 258 |
| 250 | 2 | 53 | 216024 | 135 | 6045 | 156 | 170 | 264 |
| 251 | 2 | 54 | 216024 | 136 | 6180 | 157 | 170 | 265 |
| 252 | 2 | 54 | 216060 | 137 | 6231 | 149 | 170 | 265 |
| 253 | 2 | 54 | 243396 | 138 | 6330 | 180 | 170 | 265 |
| 254 | 2 | 54 | 253072 | 139 | 6388 | 180 | 170 | 266 |
| 255 | 2 | 54 | 255318 | 140 | 6445 | 191 | 170 | 266 |
| 256 | 2 | 54 | 256868 | 136 | 6224 | 218 | 240 | 270 |
| 257 | 2 | 55 | 256868 | 137 | 6360 | 219 | 240 | 271 |
| 258 | 2 | 55 | 258868 | 138 | 6408 | 219 | 240 | 272 |
| 259 | 2 | 55 | 283609 | 139 | 6490 | 250 | 240 | 273 |
| 260 | 2 | 55 | 283912 | 140 | 6558 | 232 | 240 | 271 |
| 261 | 2 | 55 | 290690 | 139 | 6521 | 271 | 240 | 275 |
| 262 | 2 | 55 | 301830 | 140 | 6581 | 271 | 240 | 276 |
| 263 | 2 | 56 | 301830 | 141 | 6721 | 272 | 240 | 277 |
| 264 | 2 | 56 | 301871 | 142 | 6782 | 276 | 170 | 273 |
| 265 | 2 | 56 | 322260 | 143 | 6858 | 276 | 170 | 273 |
| 266 | 2 | 56 | 323777 | 144 | 6901 | 249 | 170 | 275 |
| 267 | 2 | 56 | 337684 | 145 | 6965 | 249 | 170 | 276 |
| 268 | 2 | 56 | 345360 | 145 | 6968 | 249 | 170 | 278 |
| 269 | 2 | 57 | 345360 | 146 | 7113 | 250 | 170 | 279 |
| 270 | 2 | 57 | 345402 | 147 | 7173 | 249 | 144 | 275 |
| 271 | 2 | 58 | 345402 | 148 | 7320 | 250 | 144 | 276 |
| 272 | 2 | 58 | 346933 | 147 | 7277 | 257 | 144 | 277 |
| 273 | 2 | 58 | 348858 | 148 | 7362 | 289 | 144 | 279 |
| 274 | 2 | 58 | 361776 | 149 | 7426 | 289 | 144 | 280 |
| 275 | 2 | 58 | 367511 | 150 | 7519 | 294 | 144 | 282 |
| 276 | 2 | 58 | 367941 | 151 | 7567 | 272 | 144 | 283 |
| 277 | 2 | 59 | 367941 | 152 | 7718 | 273 | 144 | 284 |
| 278 | 2 | 59 | 381162 | 153 | 7783 | 273 | 144 | 285 |
| 279 | 2 | 59 | 387657 | 152 | 7742 | 267 | 144 | 293 |
| 280 | 2 | 59 | 387705 | 153 | 7810 | 256 | 80 | 294 |
| 281 | 2 | 60 | 387705 | 154 | 7963 | 257 | 80 | 295 |
| 282 | 2 | 60 | 389923 | 155 | 8014 | 257 | 80 | 296 |
| 283 | 2 | 61 | 389923 | 156 | 8169 | 258 | 80 | 297 |
| 284 | 2 | 61 | 397705 | 156 | 8172 | 258 | 80 | 300 |
| 285 | 2 | 61 | 399378 | 157 | 8235 | 240 | 80 | 299 |
| 286 | 2 | 61 | 401042 | 158 | 8320 | 218 | 80 | 296 |
| 287 | 2 | 61 | 433661 | 159 | 8408 | 218 | 80 | 297 |
| 288 | 2 | 61 | 433714 | 159 | 8403 | 207 | 51 | 298 |
| 289 | 2 | 61 | 531152 | 150 | 7695 | 196 | 102 | 324 |
| 290 | 2 | 61 | 535198 | 151 | 7744 | 200 | 102 | 315 |
| 291 | 2 | 61 | 557710 | 152 | 7812 | 200 | 102 | 316 |
| 292 | 2 | 61 | 568212 | 152 | 7815 | 200 | 102 | 320 |
| 293 | 2 | 62 | 568212 | 153 | 7967 | 201 | 102 | 321 |
| 294 | 2 | 62 | 568966 | 153 | 7954 | 201 | 82 | 322 |
| 295 | 2 | 62 | 606838 | 154 | 8035 | 201 | 82 | 321 |
| 296 | 2 | 62 | 611758 | 154 | 8041 | 211 | 140 | 327 |
| 297 | 2 | 62 | 614368 | 155 | 8129 | 193 | 122 | 331 |
| 298 | 2 | 62 | 635774 | 156 | 8196 | 193 | 122 | 332 |
| 299 | 2 | 62 | 681858 | 157 | 8312 | 221 | 122 | 333 |
| 300 | 2 | 62 | 681936 | 157 | 8305 | 212 | 86 | 335 |
| 301 | 2 | 62 | 730048 | 158 | 8395 | 212 | 86 | 336 |
| 302 | 2 | 62 | 754948 | 159 | 8463 | 212 | 86 | 337 |
| 303 | 2 | 62 | 786658 | 160 | 8535 | 212 | 86 | 338 |
| 304 | 2 | 62 | 790518 | 159 | 8489 | 219 | 134 | 341 |
| 305 | 2 | 62 | 838930 | 160 | 8573 | 219 | 134 | 340 |
| 306 | 2 | 62 | 839111 | 161 | 8647 | 224 | 86 | 321 |
| 307 | 2 | 63 | 839111 | 162 | 8808 | 225 | 86 | 322 |
| 308 | 2 | 63 | 840733 | 163 | 8886 | 225 | 86 | 325 |
| 309 | 2 | 63 | 873931 | 164 | 8959 | 225 | 86 | 326 |
| 310 | 2 | 63 | 878399 | 165 | 9011 | 208 | 86 | 321 |
| 311 | 2 | 64 | 878399 | 166 | 9176 | 209 | 86 | 322 |
| 312 | 2 | 64 | 878974 | 166 | 9175 | 200 | 86 | 323 |
| 313 | 2 | 65 | 878974 | 167 | 9341 | 201 | 86 | 324 |
| 314 | 2 | 65 | 908988 | 168 | 9413 | 201 | 86 | 325 |
| 315 | 2 | 65 | 909441 | 169 | 9492 | 192 | 106 | 325 |
| 316 | 2 | 65 | 925219 | 169 | 9494 | 192 | 106 | 327 |
| 317 | 2 | 66 | 925219 | 170 | 9663 | 193 | 106 | 328 |
| 318 | 2 | 66 | 932411 | 171 | 9722 | 193 | 106 | 329 |
| 319 | 2 | 66 | 1076556 | 172 | 9837 | 220 | 106 | 329 |
| 320 | 2 | 66 | 1077192 | 172 | 9831 | 215 | 150 | 330 |
| 321 | 2 | 66 | 1118852 | 173 | 9907 | 215 | 150 | 331 |
| 322 | 2 | 66 | 1122887 | 174 | 9959 | 195 | 150 | 333 |
| 323 | 2 | 66 | 1250895 | 175 | 10111 | 242 | 150 | 332 |
| 324 | 2 | 66 | 1251276 | 174 | 10028 | 237 | 150 | 333 |
| 325 | 2 | 66 | 1269654 | 174 | 10059 | 262 | 150 | 337 |
| 326 | 2 | 66 | 1312203 | 175 | 10134 | 262 | 150 | 338 |
| 327 | 2 | 66 | 1359913 | 176 | 10210 | 262 | 150 | 339 |
| 328 | 2 | 66 | 1370001 | 176 | 10215 | 262 | 150 | 346 |
| 329 | 2 | 66 | 1462895 | 177 | 10310 | 262 | 150 | 347 |
| 330 | 2 | 66 | 1463110 | 178 | 10379 | 252 | 66 | 342 |
| 331 | 2 | 67 | 1463110 | 179 | 10557 | 253 | 66 | 343 |
| 332 | 2 | 67 | 1488164 | 179 | 10559 | 253 | 66 | 344 |
| 333 | 2 | 67 | 1518121 | 178 | 10513 | 293 | 66 | 351 |
| 334 | 2 | 67 | 1567446 | 179 | 10591 | 293 | 66 | 352 |
| 335 | 2 | 67 | 1665137 | 180 | 10683 | 293 | 66 | 352 |
| 336 | 2 | 67 | 1665305 | 181 | 10750 | 283 | 66 | 353 |
| 337 | 2 | 68 | 1665305 | 182 | 10931 | 284 | 66 | 354 |
| 338 | 2 | 68 | 1683042 | 181 | 10812 | 319 | 100 | 360 |
| 339 | 2 | 68 | 1742485 | 182 | 10891 | 319 | 100 | 361 |
| 340 | 2 | 68 | 1743248 | 183 | 10983 | 309 | 100 | 350 |
| 341 | 2 | 68 | 1954475 | 184 | 11103 | 309 | 100 | 350 |
| 342 | 2 | 68 | 1957398 | 185 | 11170 | 271 | 100 | 351 |
| 343 | 2 | 68 | 2064798 | 180 | 10813 | 330 | 168 | 369 |
| 344 | 2 | 68 | 2079910 | 180 | 10815 | 330 | 168 | 371 |
| 345 | 2 | 68 | 2095448 | 181 | 10878 | 344 | 168 | 381 |
| 346 | 2 | 68 | 2163864 | 182 | 10957 | 344 | 168 | 382 |
| 347 | 2 | 69 | 2163864 | 183 | 11139 | 345 | 168 | 383 |
| 348 | 2 | 69 | 2167539 | 184 | 11203 | 349 | 168 | 384 |
| 349 | 2 | 70 | 2167539 | 185 | 11387 | 350 | 168 | 385 |
| 350 | 2 | 70 | 2168306 | 186 | 11461 | 327 | 116 | 385 |
| 351 | 2 | 70 | 2176926 | 186 | 11477 | 354 | 150 | 393 |
| 352 | 2 | 70 | 2180627 | 186 | 11473 | 388 | 258 | 396 |
| 353 | 2 | 71 | 2180627 | 187 | 11659 | 389 | 258 | 397 |
| 354 | 2 | 71 | 2195316 | 188 | 11727 | 389 | 258 | 399 |
| 355 | 2 | 71 | 2326595 | 189 | 11824 | 389 | 258 | 397 |
| 356 | 2 | 71 | 2368452 | 189 | 11826 | 389 | 258 | 401 |
| 357 | 2 | 71 | 2370676 | 190 | 11939 | 379 | 258 | 392 |
| 358 | 2 | 71 | 2444745 | 191 | 12021 | 379 | 258 | 393 |
| 359 | 2 | 72 | 2444745 | 192 | 12212 | 380 | 258 | 394 |
| 360 | 2 | 72 | 2444775 | 192 | 12208 | 371 | 258 | 395 |
| 361 | 2 | 72 | 2744232 | 181 | 11132 | 283 | 258 | 433 |
| 362 | 2 | 72 | 2827980 | 182 | 11211 | 283 | 258 | 434 |
| 363 | 2 | 72 | 2858344 | 181 | 11103 | 306 | 354 | 495 |
| 364 | 2 | 72 | 2859438 | 182 | 11198 | 298 | 270 | 469 |
| 365 | 2 | 72 | 3008106 | 183 | 11290 | 298 | 270 | 465 |
| 366 | 2 | 72 | 3028916 | 184 | 11356 | 298 | 270 | 469 |
| 367 | 2 | 73 | 3028916 | 185 | 11540 | 299 | 270 | 470 |
| 368 | 2 | 73 | 3039736 | 184 | 11480 | 311 | 270 | 473 |
| 369 | 2 | 73 | 3090212 | 183 | 11434 | 311 | 270 | 488 |
| 370 | 2 | 73 | 3108982 | 184 | 11494 | 325 | 270 | 494 |
| 371 | 2 | 73 | 3295284 | 185 | 11602 | 325 | 270 | 489 |
| 372 | 2 | 73 | 3298998 | 186 | 11664 | 297 | 270 | 492 |
| 373 | 2 | 74 | 3298998 | 187 | 11850 | 298 | 270 | 493 |
| 374 | 2 | 74 | 3309930 | 188 | 11957 | 310 | 270 | 477 |
| 375 | 2 | 74 | 3325324 | 187 | 11857 | 302 | 354 | 495 |
| 376 | 2 | 74 | 3342872 | 187 | 11857 | 302 | 354 | 501 |
| 377 | 2 | 74 | 3682908 | 188 | 11998 | 363 | 354 | 501 |
| 378 | 2 | 74 | 3683946 | 189 | 12081 | 366 | 558 | 500 |
| 379 | 2 | 75 | 3683946 | 190 | 12270 | 367 | 558 | 501 |
| 380 | 2 | 75 | 3684467 | 191 | 12369 | 382 | 279 | 478 |
| 381 | 2 | 75 | 3813769 | 192 | 12458 | 382 | 279 | 479 |
| 382 | 2 | 75 | 3929399 | 193 | 12541 | 382 | 279 | 480 |
| 383 | 2 | 76 | 3929399 | 194 | 12734 | 383 | 279 | 481 |
| 384 | 2 | 76 | 3930946 | 193 | 12640 | 413 | 360 | 482 |
| 385 | 2 | 76 | 3939355 | 194 | 12763 | 433 | 258 | 474 |
| 386 | 2 | 76 | 4055259 | 195 | 12847 | 433 | 258 | 475 |
| 387 | 2 | 76 | 4113916 | 194 | 12795 | 433 | 258 | 489 |
| 388 | 2 | 76 | 4181844 | 194 | 12795 | 433 | 258 | 500 |
| 389 | 2 | 77 | 4181844 | 195 | 12989 | 434 | 258 | 501 |
| 390 | 2 | 77 | 4181999 | 196 | 13073 | 411 | 206 | 493 |
| 391 | 2 | 77 | 4715061 | 197 | 13240 | 556 | 206 | 494 |
| 392 | 2 | 77 | 4718688 | 197 | 13225 | 521 | 162 | 499 |
| 393 | 2 | 77 | 4876591 | 198 | 13316 | 521 | 162 | 500 |
| 394 | 2 | 77 | 5017175 | 199 | 13403 | 521 | 162 | 501 |
| 395 | 2 | 77 | 5276525 | 200 | 13509 | 521 | 162 | 497 |
| 396 | 2 | 77 | 5276873 | 201 | 13594 | 507 | 162 | 494 |
| 397 | 2 | 78 | 5276873 | 202 | 13795 | 508 | 162 | 495 |
| 398 | 2 | 78 | 5426777 | 203 | 13883 | 508 | 162 | 496 |
| 399 | 2 | 78 | 5431954 | 204 | 14010 | 551 | 129 | 463 |
| 400 | 2 | 78 | 5432967 | 202 | 13806 | 520 | 96 | 466 |
| 401 | 2 | 79 | 5432967 | 203 | 14008 | 521 | 96 | 467 |
| 402 | 2 | 79 | 5465217 | 204 | 14077 | 521 | 96 | 465 |
| 403 | 2 | 79 | 5773678 | 205 | 14222 | 442 | 96 | 466 |
| 404 | 2 | 79 | 5868487 | 205 | 14226 | 442 | 96 | 475 |
| 405 | 2 | 79 | 5876541 | 206 | 14323 | 456 | 120 | 483 |
| 406 | 2 | 79 | 5906571 | 207 | 14388 | 505 | 120 | 484 |
| 407 | 2 | 79 | 6865807 | 208 | 14529 | 635 | 120 | 472 |
| 408 | 2 | 79 | 6867693 | 209 | 14618 | 596 | 120 | 473 |
| 409 | 2 | 80 | 6867693 | 210 | 14827 | 597 | 120 | 474 |
| 410 | 2 | 80 | 6906334 | 211 | 14902 | 633 | 120 | 470 |
| 411 | 2 | 80 | 7140833 | 212 | 14997 | 633 | 120 | 471 |
| 412 | 2 | 80 | 7251934 | 212 | 15001 | 633 | 120 | 478 |
| 413 | 2 | 80 | 7685552 | 213 | 15123 | 633 | 120 | 479 |
| 414 | 2 | 80 | 7697680 | 214 | 15196 | 586 | 120 | 492 |
| 415 | 2 | 80 | 8095626 | 215 | 15307 | 586 | 120 | 490 |
| 416 | 2 | 80 | 8104530 | 216 | 15409 | 565 | 120 | 491 |
| 417 | 2 | 80 | 8372689 | 217 | 15506 | 565 | 120 | 492 |
| 418 | 2 | 80 | 8394698 | 218 | 15633 | 592 | 120 | 473 |
| 419 | 2 | 81 | 8394698 | 219 | 15851 | 593 | 120 | 474 |
| 420 | 2 | 81 | 8394754 | 220 | 15934 | 581 | 110 | 472 |
| 421 | 2 | 82 | 8394754 | 221 | 16154 | 582 | 110 | 473 |
| 422 | 2 | 82 | 8629471 | 222 | 16249 | 582 | 110 | 474 |
| 423 | 2 | 82 | 8742239 | 221 | 16184 | 548 | 86 | 492 |
| 424 | 2 | 82 | 8792920 | 221 | 16184 | 582 | 110 | 494 |
| 425 | 2 | 82 | 8882557 | 222 | 16321 | 596 | 114 | 504 |
| 426 | 2 | 82 | 8938237 | 223 | 16398 | 596 | 114 | 501 |
| 427 | 2 | 82 | 9413895 | 224 | 16523 | 596 | 114 | 502 |
| 428 | 2 | 82 | 9554721 | 224 | 16526 | 596 | 114 | 514 |

| 429 | 2 | 82 | 9578133 | 225 | 16658 | 533 | 114 | 511 |
| --- | --- | --- | --- | --- | --- | --- | --- | --- |
| 430 | 2 | 82 | 9622020 | 226 | 16736 | 533 | 114 | 501 |
| 431 | 2 | 83 | 9622020 | 227 | 16962 | 534 | 114 | 502 |
| 432 | 2 | 83 | 9622744 | 226 | 16870 | 513 | 110 | 503 |
| 433 | 2 | 84 | 9622744 | 227 | 17096 | 514 | 110 | 504 |
| 434 | 2 | 84 | 9645552 | 228 | 17163 | 461 | 110 | 506 |
| 435 | 2 | 84 | 9690786 | 229 | 17252 | 493 | 110 | 506 |
| 436 | 2 | 84 | 9821579 | 229 | 17255 | 493 | 110 | 516 |
| 437 | 2 | 84 | 10561514 | 230 | 17449 | 554 | 110 | 516 |
| 438 | 2 | 84 | 10630383 | 231 | 17531 | 554 | 110 | 513 |
| 439 | 2 | 85 | 10630383 | 232 | 17762 | 555 | 110 | 514 |
| 440 | 2 | 85 | 10632891 | 233 | 17856 | 555 | 110 | 515 |
| 441 | 2 | 85 | 10644013 | 230 | 17528 | 541 | 122 | 498 |
| 442 | 2 | 85 | 10675730 | 231 | 17653 | 497 | 110 | 500 |
| 443 | 2 | 86 | 10675730 | 232 | 17884 | 498 | 110 | 501 |
| 444 | 2 | 86 | 10684060 | 233 | 17963 | 471 | 110 | 501 |
| 445 | 2 | 86 | 11205356 | 234 | 18081 | 471 | 110 | 496 |
| 446 | 2 | 86 | 11502000 | 235 | 18181 | 471 | 110 | 497 |
| 447 | 2 | 86 | 11877069 | 236 | 18290 | 471 | 110 | 498 |
| 448 | 2 | 86 | 11880755 | 236 | 18283 | 470 | 110 | 499 |

# Appendix B: F Vactors of the polytopes representing the natural numbers from 1 to 69

| N | f₀ | f₁ | f₂ | f₃ | f₄ | f₅ | f₆ | f₇ | f₈ | f₉ | f₁₀ | f₁₁ | f₁₂ | f₁₃ | f₁₄ | f₁₅ | f₁₆ | f₁₇ | f₁₈ |
|---|----|----|----|----|----|----|----|----|----|----|-----|-----|-----|-----|-----|-----|-----|-----|-----|
| 1 | | | | | | | | | | | | | | | | | | | |
| 2 | 2 | | | | | | | | | | | | | | | | | | |
| 3 | 3 | 3 | | | | | | | | | | | | | | | | | |
| 4 | 3 | 3 | | | | | | | | | | | | | | | | | |
| 5 | 4 | 6 | 4 | | | | | | | | | | | | | | | | |
| 6 | 5 | 8 | 5 | | | | | | | | | | | | | | | | |
| 7 | 6 | 13 | 13 | 6 | | | | | | | | | | | | | | | |
| 8 | 6 | 13 | 13 | 6 | | | | | | | | | | | | | | | |
| 9 | 5 | 10 | 10 | 5 | | | | | | | | | | | | | | | |
| 10 | 6 | 13 | 13 | 6 | | | | | | | | | | | | | | | |
| 11 | 7 | 19 | 26 | 19 | 7 | | | | | | | | | | | | | | |
| 12 | 8 | 23 | 32 | 23 | 8 | | | | | | | | | | | | | | |
| 13 | 9 | 31 | 55 | 55 | 31 | 9 | | | | | | | | | | | | | |
| 14 | 10 | 35 | 61 | 59 | 32 | 9 | | | | | | | | | | | | | |
| 15 | 11 | 40 | 71 | 69 | 37 | 10 | | | | | | | | | | | | | |
| 16 | 10 | 36 | 65 | 65 | 36 | 10 | | | | | | | | | | | | | |
| 17 | 11 | 46 | 101 | 130 | 101 | 46 | 11 | | | | | | | | | | | | |
| 18 | 12 | 51 | 111 | 140 | 106 | 47 | 11 | | | | | | | | | | | | |
| 19 | 13 | 63 | 162 | 251 | 246 | 153 | 58 | 12 | | | | | | | | | | | |
| 20 | 13 | 63 | 162 | 251 | 246 | 153 | 58 | 12 | | | | | | | | | | | |
| 21 | 14 | 70 | 183 | 286 | 281 | 174 | 65 | 13 | | | | | | | | | | | |
| 22 | 15 | 76 | 198 | 306 | 296 | 180 | 66 | 13 | | | | | | | | | | | |
| 23 | 16 | 91 | 274 | 504 | 602 | 476 | 246 | 79 | 14 | | | | | | | | | | |
| 24 | 17 | 97 | 288 | 518 | 602 | 462 | 232 | 73 | 13 | | | | | | | | | | |
| 25 | 15 | 83 | 246 | 448 | 532 | 420 | 218 | 71 | 13 | | | | | | | | | | |
| 26 | 16 | 89 | 261 | 468 | 547 | 426 | 219 | 71 | 13 | | | | | | | | | | |
| 27 | 15 | 83 | 246 | 448 | 532 | 420 | 218 | 71 | 13 | | | | | | | | | | |
| 28 | 15 | 84 | 253 | 469 | 567 | 455 | 239 | 78 | 14 | | | | | | | | | | |
| 29 | 16 | 99 | 337 | 722 | 1036 | 1022 | 694 | 317 | 92 | 15 | | | | | | | | | |
| 30 | 17 | 108 | 372 | 799 | 1141 | 1113 | 743 | 332 | 94 | 15 | | | | | | | | | |
| 31 | 18 | 125 | 480 | 1171 | 1940 | 2254 | 1856 | 1075 | 426 | 109 | 16 | | | | | | | | |
| 32 | 17 | 116 | 444 | 1087 | 1814 | 2128 | 1772 | 1039 | 417 | 108 | 16 | | | | | | | | |
| 33 | 18 | 126 | 489 | 1207 | 2024 | 2380 | 1982 | 1159 | 462 | 118 | 17 | | | | | | | | |
| 34 | 19 | 135 | 525 | 1291 | 2150 | 2506 | 2066 | 1195 | 471 | 119 | 17 | | | | | | | | |
| 35 | 20 | 149 | 606 | 1555 | 2696 | 3262 | 2780 | 1651 | 660 | 165 | 22 | | | | | | | | |
| 36 | 20 | 148 | 597 | 1519 | 2612 | 3136 | 2654 | 1567 | 624 | 156 | 21 | | | | | | | | |
| 37 | 21 | 168 | 745 | 2116 | 4131 | 5748 | 5790 | 4221 | 2191 | 780 | 177 | 22 | | | | | | | |
| 38 | 22 | 178 | 789 | 2228 | 4313 | 5944 | 5930 | 4285 | 2208 | 782 | 177 | 22 | | | | | | | |
| 39 | 23 | 189 | 844 | 2392 | 4635 | 6378 | 6336 | 4545 | 2317 | 809 | 180 | 22 | | | | | | | |
| 40 | 24 | 200 | 896 | 2531 | 4865 | 6616 | 6476 | 4567 | 2286 | 784 | 172 | 21 | | | | | | | |
| 41 | 25 | 224 | 1096 | 3427 | 7396 | 11481 | 13092 | 11043 | 6853 | 3070 | 956 | 193 | 22 | | | | | | |
| 42 | 26 | 233 | 1130 | 3492 | 7441 | 11403 | 12840 | 10701 | 6568 | 2915 | 902 | 182 | 21 | | | | | | |
| 43 | 27 | 259 | 1363 | 4622 | 10933 | 18844 | 24243 | 23541 | 17269 | 9483 | 3817 | 1084 | 203 | 22 | | | | | |
| 44 | 27 | 259 | 1364 | 4632 | 10978 | 18964 | 24453 | 23793 | 17479 | 9603 | 3862 | 1094 | 204 | 22 | | | | | |
| 45 | 28 | 275 | 1477 | 5106 | 12308 | 21616 | 28335 | 28029 | 20932 | 11683 | 4763 | 1360 | 252 | 26 | | | | | |
| 46 | 29 | 287 | 1541 | 5307 | 12722 | 22204 | 28923 | 28443 | 21133 | 11747 | 4775 | 1361 | 252 | 26 | | | | | |
| 47 | 30 | 316 | 1828 | 6848 | 18029 | 34926 | 51127 | 57366 | 49576 | 32880 | 16522 | 6136 | 1613 | 278 | 27 | | | | |
| 48 | 30 | 314 | 1803 | 6704 | 17523 | 33716 | 49048 | 54726 | 47068 | 31098 | 15587 | 5784 | 1523 | 264 | 26 | | | | |
| 49 | 26 | 259 | 1455 | 5359 | 13971 | 26928 | 39358 | 44220 | 38368 | 25619 | 13003 | 4899 | 1315 | 234 | 24 | | | | |
| 50 | 26 | 255 | 1406 | 5082 | 13014 | 24673 | 35530 | 39402 | 33814 | 22385 | 11298 | 4250 | 1146 | 207 | 22 | | | | |
| 51 | 27 | 267 | 1472 | 5302 | 13509 | 25465 | 36454 | 40194 | 34309 | 22605 | 11364 | 4262 | 1147 | 207 | 22 | | | | |
| 52 | 27 | 268 | 1484 | 5368 | 13729 | 25960 | 37246 | 41118 | 35101 | 23100 | 11584 | 4328 | 1159 | 208 | 22 | | | | |
| 53 | 28 | 295 | 1752 | 6852 | 19097 | 39689 | 63206 | 78364 | 76219 | 58201 | 34684 | 15912 | 5487 | 1367 | 230 | 23 | | | |
| 54 | 29 | 308 | 1828 | 7113 | 19668 | 40470 | 63712 | 78001 | 74866 | 56408 | 33188 | 15055 | 5148 | 1278 | 216 | 22 | | | |
| 55 | 30 | 329 | 2013 | 8062 | 22905 | 48335 | 77869 | 97306 | 95029 | 72567 | 43055 | 19566 | 6643 | 1617 | 263 | 25 | | | |
| 56 | 31 | 346 | 2147 | 8712 | 25076 | 53626 | 87593 | 111034 | 110044 | 85294 | 51349 | 23648 | 8112 | 1982 | 319 | 29 | | | |
| 57 | 32 | 360 | 2239 | 9088 | 26143 | 55848 | 91091 | 115258 | 113971 | 88088 | 52845 | 24232 | 8269 | 2008 | 321 | 29 | | | |
| 58 | 33 | 374 | 2329 | 9440 | 27078 | 57630 | 93599 | 117898 | 116050 | 89298 | 53351 | 24376 | 8294 | 2010 | 321 | 29 | | | |
| 59 | 34 | 407 | 2703 | 11769 | 36518 | 84708 | 151229 | 211497 | 233948 | 205348 | 142649 | 77727 | 32670 | 10304 | 2331 | 350 | 30 | | |
| 60 | 35 | 422 | 2808 | 12224 | 37883 | 87711 | 156234 | 217932 | 240383 | 210353 | 145652 | 79092 | 33125 | 10409 | 2346 | 351 | 30 | | |
| 61 | 36 | 457 | 3230 | 15032 | 50107 | 125594 | 243945 | 374166 | 458315 | 450736 | 356005 | 224744 | 112217 | 43534 | 12755 | 2697 | 381 | 31 | |
| 62 | 37 | 472 | 3334 | 15474 | 51394 | 128311 | 248235 | 379314 | 463034 | 454025 | 357721 | 225394 | 112386 | 43561 | 12757 | 2697 | 381 | 31 | |
| 63 | 38 | 494 | 3545 | 16680 | 56063 | 141415 | 275990 | 424788 | 521521 | 513513 | 405626 | 255788 | 127401 | 49217 | 14326 | 2999 | 417 | 33 | |
| 64 | 36 | 461 | 3290 | 15455 | 51968 | 131314 | 256971 | 396903 | 489346 | 484198 | 384605 | 244049 | 122396 | 47642 | 13981 | 2952 | 414 | 33 | |
| 65 | 37 | 483 | 3495 | 16575 | 56062 | 142039 | 277914 | 428077 | 525096 | 515801 | 405912 | 254709 | 126114 | 48397 | 13990 | 2910 | 403 | 32 | |
| 66 | 38 | 495 | 3556 | 16729 | 56153 | 141311 | 274911 | 421499 | 515229 | 504933 | 396903 | 249067 | 123475 | 47501 | 13781 | 2880 | 401 | 32 | |
| 67 | 39 | 533 | 4051 | 20285 | 72882 | 197464 | 416222 | 696410 | 936728 | 1020162 | 901836 | 645970 | 372542 | 170976 | 61282 | 16661 | 3281 | 433 | 33 |
| 68 | 39 | 535 | 4086 | 20572 | 74347 | 202679 | 429963 | 724165 | 980629 | 1075217 | 956891 | 689871 | 400297 | 184717 | 66497 | 18126 | 3568 | 468 | 35 |
| 69 | 40 | 552 | 4224 | 21280 | 76909 | 209595 | 444341 | 747617 | 1010945 | 1106391 | 982345 | 706251 | 408487 | 187825 | 67359 | 18290 | 3587 | 469 | 35 |